# Invariance principles for fractionally integrated nonlinear processes[*]

## Wei Biao Wu[1] and Xiaofeng Shao[2]

*University of Chicago and University Illinois, at Urbana-Champaign*

**Abstract:** We obtain invariance principles for a wide class of fractionally integrated nonlinear processes. The limiting distributions are shown to be fractional Brownian motions. Under very mild conditions, we extend earlier ones on long memory linear processes to a more general setting. The invariance principles are applied to the popular R/S and KPSS tests.

## 1. Introduction

Invariance principles (or functional central limit theorems) play an important role in econometrics and statistics. For example, to obtain asymptotic distributions of unit-root test statistics, researchers have applied invariance principles of various forms; see [24, 30, 40] among others. The primary goal of this paper is to establish invariance principles for a class of fractionally integrated nonlinear processes. Let the process

$$(1.1) \qquad u_t = F(\ldots, \varepsilon_{t-1}, \varepsilon_t), \quad t \in \mathbf{Z},$$

where $\varepsilon_t$ are independent and identically distributed (iid) random variables and $F$ is a measurable function such that $u_t$ is well-defined. Then $u_t$ is stationary and causal. Let $d \in (-1/2, 1/2)$ and define the Type I fractional I($d$) process $X_t$ by

$$(1.2) \qquad (1-B)^d (X_t - \mu) = u_t, \ t \in \mathbf{Z},$$

where $\mu$ is the mean and $B$ is the backward shift operator: $BX_t = X_{t-1}$. The Type II I($d$) fractional process is defined as

$$(1.3) \qquad (1-B)^d (Y_t - Y_0) = u_t \mathbf{1}(t \geq 1).$$

where $Y_0$ is a random variable whose distribution is independent of $t$. There has been a recent surge of interest in Type II processes [25, 28] and it arises naturally when the processes start at a given time point. The framework (1.1) includes a very wide class of processes [26, 28, 33, 34, 37]. It includes linear processes $u_t = \sum_{j=0}^{\infty} b_j \varepsilon_{t-j}$ as a special case. It also includes a large class of nonlinear time series models, such as bilinear models, threshold models and GARCH type models [29, 41]. Recently, fractionally integrated autoregressive and moving average models (FARIMA) with

[*]The work is supported in part by NSF Grant DMS-0478704.
[1]Department of Statistics, University of Chicago, 5734 S. University Ave, Chicago, IL 60637, e-mail: `wbwu@galton.uchicago.edu`
[2]Department of Statistics, University of Illinois, at Urbana-Champaign, 725 Wright St, Champaign, IL 61820, e-mail: `xshao@uiuc.edu`

GARCH innovations have attracted much attention in financial time series modeling (see [2]). In financial time series analysis, conditional heteroscedasticity and long memory are commonly seen [14, 19]. The FARIMA-GARCH model naturally fits into our framework.

Most of the results in the literature assume $\{u_t\}$ to be either iid or linear processes. Recently, Wu and Min [41] established an invariance principle under (1.2) when $d \in [0, 1/2)$. The literature seems more concentrated on the case $d \in (0, 1/2)$. Part of the reason is that this case corresponds to long memory and it appears in various areas such as finance, hydrology and telecommunication. When $d \in (1/2, 1)$, the process is non-stationary and it can be defined as $\sum_{s=1}^{t} X_s$ or $\sum_{s=1}^{t} Y_s$, where $X_s$ and $Y_s$ are Type I and Type II I$(d-1)$ processes, respectively. Empirical evidence of $d \in (1/2, 1)$ has been found by Byers et al. [4] in poll data modeling and Kim [16] in macroeconomic time series. Therefore the study of partial sums of I$(d)$, $d \in (-1/2, 0)$ is also of interest since it naturally leads to I$(d)$ processes, $d \in (1/2, 1)$. In fact, our results can be easily extended to the process with fractional index $p + d$, $p \in \mathbf{N}$, $d \in (-1/2, 0) \cup (0, 1/2)$ [cf. Corollary 2.1].

The study of the invariance principle has a long history. Here we only mention some representatives: Davydov [8], Mcleich [23], Gorodetskii [11], Hall and Heyde [12], Davidson and De Jong [6], De Jong and Davidson [7] and the references cited therein. Most of them deal with Type I processes. Recent developments for Type II processes can be found in [15, 22, 36] among others.

The paper is structured as follows. Section 2 presents invariance principles for both types of processes. Section 3 considers limit distributions of tests of long memory under mild moment conditions. Technical details are given in the appendix.

## 2. Main Results

We first define two types of fractional Brownian motions. For Type I fractional Brownian motion, let $d \in (-1/2, 1/2)$ and

$$\mathbb{B}_d(t) = \frac{1}{A(d)} \int_{-\infty}^{\infty} [\{(t-s)_+\}^d - \{(-s)_+\}^d] d\mathbb{B}(s), \quad t \in \mathbf{R},$$

where $(t)_+ = \max(t, 0)$, $\mathbb{B}(s)$ is a standard Brownian motion and

$$A(d) = \left\{ \frac{1}{2d+1} + \int_0^{\infty} [(1+s)^d - s^d]^2 ds \right\}^{1/2}.$$

Type II fractional Brownian motion $\{W_d(t), t \geq 0\}$, $d > -1/2$, is defined as

$$W_d(0) = 0, \quad W_d(t) = (2d+1)^{1/2} \int_0^t (t-s)^d d\mathbb{B}(s).$$

The main difference of $\mathbb{B}_d(t)$ and $W_d(t)$ lies in the prehistoric treatment. See [22] for a detailed discussion of the difference between them. Here we are interested in the weak convergence of the partial sums $T_m = \sum_{i=1}^m X_i$ and $T_m^\diamond = \sum_{i=1}^m Y_i$. Let $\mathcal{D}[0, 1]$ be the space of functions on $[0, 1]$ which are right continuous and have left-hand limits, endowed with the Skorohod topology [3]. Denote weak convergence by "$\Rightarrow$".

For a random variable $X$, write $X \in \mathcal{L}^p$ ($p > 0$) if $\|X\|_p := [\mathbf{E}(|X|^p)]^{1/p} < \infty$ and $\|\cdot\| = \|\cdot\|_2$. Let $\mathcal{F}_t = (\ldots, \varepsilon_{t-1}, \varepsilon_t)$ be the shift process. Define the projections $\mathcal{P}_k$ by $\mathcal{P}_k X = \mathbf{E}(X|\mathcal{F}_k) - \mathbf{E}(X|\mathcal{F}_{k-1})$, $X \in \mathcal{L}^1$. For two sequences $(a_n)$, $(b_n)$, denote by $a_n \sim b_n$ if $a_n/b_n \to 1$ as $n \to \infty$. The symbols "$\to_D$" and "$\to_\mathbf{P}$" stand for convergence in distribution and in probability, respectively. The symbols $O_\mathbf{P}(1)$ and $o_\mathbf{P}(1)$ signify being bounded in probability and convergence to zero in probability. Let $N(\mu, \sigma^2)$ be a normal distribution with mean $\mu$ and variance $\sigma^2$. Hereafter we assume without loss of generality that $\mathbf{E}(u_t) = 0$, $\mu = 0$ and $Y_0 = 0$. Let $\{\varepsilon'_k, k \in \mathbf{Z}\}$ be an iid copy of $\{\varepsilon_k, k \in \mathbf{Z}\}$ and $\mathcal{F}^*_k = (\mathcal{F}_{-1}, \varepsilon'_0, \varepsilon_1, \ldots, \varepsilon_k)$.

Theorems 2.1 and 2.2 concern Type I and II processes respectively. Using the continuous mapping theorem, Theorems 2.1 and 2.2 imply Corollary 2.1 which deals with general fractional processes with higher orders. For $d \in (-1/2, 0)$, an undesirable feature of our results is that the moment condition depends on $d$. However, this seems to be necessary; see Remark 4.1. Similar conditions were imposed in [30, 36]. Theorem 2.2 extends early results by Akonom and Gourieroux [1], Tanaka [31] and Wang et al. [36], who assumed $u_t$ to be either iid or linear processes. See [15, 22] for a multivariate extension.

**Theorem 2.1.** *Let $u_t \in \mathcal{L}^q$, where $q > 2/(2d+1)$ if $d \in (-1/2, 0)$ and $q = 2$ if $d \in (0, 1/2)$. Assume*

$$\sum_{k=0}^{\infty} \|\mathcal{P}_0 u_k\|_q < \infty. \tag{2.1}$$

*Then $\zeta_j := \sum_{k=j}^{\infty} \mathcal{P}_j u_k \in \mathcal{L}^q$ and, if $\|\zeta_0\| > 0$,*

$$\frac{T_{\lfloor nt \rfloor}}{n^{d+1/2}} \Rightarrow \kappa_1(d)\mathbb{B}_d(t) \text{ in } \mathcal{D}[0,1], \text{ where } \kappa_1(d) = \frac{A(d)\|\zeta_0\|}{\Gamma(d+1)}. \tag{2.2}$$

**Remark 2.1.** Note that $\|\zeta_0\|^2 = 2\pi f_u(0)$, where $f_u(\cdot)$ is the spectral density function of $\{u_t\}$; see Wu [39, 41] for the details.

**Theorem 2.2.** *Under the conditions of Theorem 2.1, for Type II processes, we have*

$$\frac{T^\diamond_{\lfloor nt \rfloor}}{n^{d+1/2}} \Rightarrow \kappa_2(d)W_d(t) \text{ in } \mathcal{D}[0,1], \text{ where } \kappa_2(d) = \frac{\|\zeta_0\|(2d+1)^{-1/2}}{\Gamma(d+1)}. \tag{2.3}$$

By the continuous mapping theorem and the standard argument in [36], we have the following corollary.

**Corollary 2.1.** *Let $u_t$ satisfy conditions in Theorem 2.1; let $d \in (0, 1/2) \cup (-1/2, 0)$ and $p \in \mathbf{N}$. [a] Define the process $\tilde{X}_t$ by $(1-B)^{p+d}\tilde{X}_t = u_t$. Then*

(i) $n^{-(d+p-1/2)}\tilde{X}_{\lfloor nt \rfloor} \Rightarrow \kappa_1(d)\mathbb{B}_{d,p}(t)$ *in* $\mathcal{D}[0,1]$;
(ii) $n^{-(d+p+1/2)}\sum_{j=1}^{\lfloor nt \rfloor} \tilde{X}_j \Rightarrow \kappa_1(d)\mathbb{B}_{d,p+1}(t)$ *in* $\mathcal{D}[0,1]$;
(iii) $n^{-2(d+p)}\sum_{j=1}^{\lfloor nt \rfloor} \tilde{X}_j^2 \Rightarrow \kappa_1^2(d)\int_0^t [\mathbb{B}_{d,p}(s)]^2 ds$ *in* $\mathcal{D}[0,1]$.

*Here $\mathbb{B}_{d,p}(t)$ is defined as*

$$\mathbb{B}_{d,p}(t) = \begin{cases} \mathbb{B}_d(t), & p = 1, \\ \int_0^t \int_0^{t_{p-1}} \cdots \int_0^{t_2} \mathbb{B}_d(t_1) dt_1 dt_2 \cdots dt_{p-1}, & p \geq 2. \end{cases}$$

*[b] Define the process $\tilde{Y}_t$ by $(1-B)^{p+d}\tilde{Y}_t = u_t \mathbf{1}(t \geq 1)$. Then similarly*

(i) $n^{-(d+p-1/2)}\tilde{Y}_{\lfloor nt \rfloor} \Rightarrow \kappa_2(d)W_{d+p-1}(t)$ *in* $\mathcal{D}[0,1]$;

(ii) $n^{-(d+p+1/2)} \sum_{j=1}^{\lfloor nt \rfloor} \tilde{Y}_j \Rightarrow \kappa_2(d) W_{d+p}(t)$ *in* $\mathcal{D}[0,1]$;

(iii) $n^{-2(d+p)} \sum_{j=1}^{\lfloor nt \rfloor} \tilde{Y}_j^2 \Rightarrow \kappa_2^2(d) \int_0^t [W_{d+p-1}(s)]^2 ds$ *in* $\mathcal{D}[0,1]$.

We now discuss condition (2.1). Let $g_k(\mathcal{F}_0) = \mathbf{E}(u_k|\mathcal{F}_0)$ and $\delta_q(k) = \|g_k(\mathcal{F}_0) - g_k(\mathcal{F}_0^*)\|_q$. Then $\delta_q(k)$ measures the contribution of $\varepsilon_0$ in predicting $u_k$. In [38] it is called the *predictive dependence measure*. Since $\|\mathcal{P}_0 u_k\|_q \leq \delta_q(k) \leq 2\|\mathcal{P}_0 u_k\|_q$ [38], (2.1) is equivalent to the $q$-stability condition [38] $\sum_{k=0}^{\infty} \delta_q(k) < \infty$, which suggests short-range dependence in that the cumulative contribution of $\varepsilon_0$ in predicting future values of $u_k$ is finite. For a variety of nonlinear time series, $\delta_q(k) = O(\rho^k)$ for some $\rho \in (0,1)$. The latter is equivalent to the *geometric moment contraction (GMC)* [41, 42]. Shao and Wu [29] verified GMC for GARCH$(r,s)$ model and its asymmetric variants and showed that the GMC property is preserved under an ARMA filter. In the special case $u_t = \sum_{k=0}^{\infty} b_k \varepsilon_{t-k}$, (2.1) holds if $\sum_{k=0}^{\infty} |b_k| < \infty$ and $\varepsilon_1 \in \mathcal{L}^q$.

## 3. Applications

There has been much work on the test of long memory under the short memory null hypothesis, i.e. I(0) versus I($d$), $d \in (0, 1/2)$. For example, Lo [20] introduced modified R/S test statistics, which admits the following form:

$$Q_n = \frac{1}{w_{n,l}} \left\{ \max_{1 \leq k \leq n} \sum_{j=1}^k (X_j - \bar{X}_n) - \min_{1 \leq k \leq n} \sum_{j=1}^k (X_j - \bar{X}_n) \right\},$$

where $\bar{X}_n = n^{-1} \sum_{j=1}^n X_j$ is the sample mean and $w_{n,l}$ is the long run variance estimator of $X_t$. Following Lo [20],

$$(3.1) \qquad w_{n,l}^2 = \frac{1}{n} \sum_{j=1}^n (X_j - \bar{X}_n)^2 + 2 \sum_{j=1}^l \left(1 - \frac{j}{l+1}\right) \hat{\gamma}_j,$$

where $\hat{\gamma}_j = \frac{1}{n} \sum_{i=1}^{n-j} (X_i - \bar{X}_n)(X_{i+j} - \bar{X}_n)$, $0 \leq j < n$. The form (3.1) is equivalent to the nonparametric spectral density estimator of $\{X_t\}$ evaluated at zero frequency with Bartlett window (up to a constant factor). Here the bandwidth satisfies

$$(3.2) \qquad l = l(n) \to \infty \text{ and } l/n \to 0, \text{ as } n \to \infty.$$

Lee and Schmidt applied [18] the KPSS test [17] for I(0) versus I($d$), $d \in (-1/2, 0) \cup (0, 1/2)$. The test statistics has the form:

$$K_n = \frac{1}{w_{n,l}^2 n^2} \sum_{k=1}^n \left(\sum_{j=1}^k (X_j - \bar{X}_n)\right)^2$$

with $w_{n,l}^2$ given by (3.1). Lee and Schmidt showed that the test is consistent against fractional alternatives and derived its asymptotic distribution under the assumption that $u_t$ are iid normal random variables. Giraitis et al. investigated [10] the theoretical performance of various forms of nonparametric tests under both short memory hypotheses and long memory alternatives. In a quite general setting, we obtain asymptotic distributions of R/S and KPSS test statistics under fractional alternatives.

**Theorem 3.1.** *Suppose that $X_t$ is generated from (1.2) and $u_t$ satisfies (2.1) with some $q > \max(2, 2/(2d+1))$. Assume (3.2). Then*

$$(3.3) \qquad l^{-2d} w_{n,l}^2 \to_{\mathbf{P}} \kappa_1^2(d).$$

*Consequently, we have*

$$(3.4) \qquad \frac{l^d}{n^{d+1/2}} Q_n \to_D \sup_{0 \le t \le 1} \tilde{\mathbb{B}}_d(t) - \inf_{0 \le t \le 1} \tilde{\mathbb{B}}_d(t),$$

*where $\tilde{\mathbb{B}}_d(t)$ is the fractional Brownian bridge, i.e. $\tilde{\mathbb{B}}_d(t) = \mathbb{B}_d(t) - t\mathbb{B}_d(1)$, and*

$$(3.5) \qquad \frac{l^{2d}}{n^{2d}} K_n \to_D \int_0^1 (\tilde{B}_d(t))^2 dt.$$

**Remark 3.1.** For $d \in (0, 1/2)$, Giraitis et al. obtained [10] (3.4) under the joint cumulant summability condition

$$(3.6) \qquad \sup_h \sum_{r,s=-n}^n |\text{cum}(X_0, X_h, X_r, X_s)| = O(n^{2d}).$$

For linear processes, (3.6) can be verified. But for nonlinear fractional processes (1.2), it seems hard to directly verify (3.6). In contrast, we only need to impose a $q$-th ($q > 2$) moment condition when $d \in (0, 1/2)$. Our dependence condition (2.1) can be easily verified for various nonlinear time series models [29, 41].

## 4. Appendix

**Lemma 4.1.** *Let $a_i = i^{-\beta} \ell(i)$, $i \ge 1$, where $\ell$ is a slowly varying function and $\beta > 1/2$; let $q > (3/2 - \beta)^{-1}$ if $1 < \beta < 3/2$ and $q = 2$ if $1/2 < \beta < 1$; let $\sigma_n = A(1-\beta) n^{3/2-\beta} \ell(n)/(1-\beta)$ and $\tilde{\sigma}_n = (3 - 2\beta)^{-1/2} n^{3/2-\beta} \ell(n)/(1-\beta)$. Assume either $(1^\circ)$ $1 < \beta < 3/2$, $\sum_{i=0}^\infty a_i = 0$ or $(2^\circ)$ $1/2 < \beta < 1$. Let*

$$(4.1) \qquad \eta_j = G(\ldots, \varepsilon_{j-1}, \varepsilon_j), \; j \in \mathbf{Z},$$

*be a martingale difference sequence with respect to $\sigma(\ldots, \varepsilon_{j-1}, \varepsilon_j)$. Assume $\eta_j \in \mathcal{L}^q$. Let $Y_j = \sum_{i=0}^\infty a_i \eta_{j-i}$, $Y_j^\diamond = \sum_{i=0}^{j-1} a_i \eta_{j-i}$, $S_i = \sum_{j=1}^i Y_j$, $S_i^\diamond = \sum_{j=1}^i Y_j^\diamond$. Then we have [a] $\sigma_n^{-1} S_{\lfloor nt \rfloor} \Rightarrow \|\eta_0\| \mathbb{B}_{3/2-\beta}(t)$ in $\mathcal{D}[0,1]$ and [b] $\tilde{\sigma}_n^{-1} S_{\lfloor nt \rfloor}^\diamond \Rightarrow \|\eta_0\| W_{3/2-\beta}(t)$ in $\mathcal{D}[0,1]$.*

*Proof of Lemma 4.1.* [a] Consider $(1^\circ)$ first. For the finite dimensional convergence, we shall apply the Cramer–Wold device. Fix $0 \le t_1 < t_2 \le 1$ and let $m_1 = \lfloor nt_1 \rfloor$ and $m_2 = \lfloor nt_2 \rfloor$. Let $A_i = \sum_{j=0}^i a_j$ if $i \ge 0$ and $A_i = 0$ if $i < 0$. For $\lambda, \mu \in \mathbf{R}$ let

$$c_{n,l} = \frac{\lambda(A_{m_1-m_2+l} - A_{l-m_2}) + \mu(A_l - A_{l-m_2})}{\sigma_n},$$
$$\sigma_{\lambda\mu}^2 = [\lambda^2 t_1^{3-2\beta} + \mu^2 t_2^{3-2\beta} + \lambda\mu(t_1^{3-2\beta} + t_2^{3-2\beta} - (t_2 - t_1)^{3-2\beta})].$$

Then $(\lambda S_{m_1} + \mu S_{m_2})/\sigma_n = \sum_{l=0}^\infty c_{n,l} \eta_{m_2-l}$ has martingale difference summands and we can apply the martingale central limit theorem. By Karamata's Theorem, $A_n = -\sum_{j=n+1}^\infty a_j \sim n^{1-\beta} \ell(n)/(\beta-1)$. Elementary calculations show that

$$(4.2) \qquad \sum_{l=0}^\infty c_{n,l}^2 \to \sigma_{\lambda\mu}^2 \text{ and } \sup_{l \ge 0} |c_{n,l}| \to 0 \text{ as } n \to \infty.$$

Let $V_l = \mathbf{E}(\eta_{m_2-l}^2|\mathcal{F}_{m_2-l-1})$. By the argument in the proof of Theorem 1 in [13], (4.2) implies that $\sum_{l=0}^{\infty} c_{n,l}^2 V_l \to \sigma_{\lambda\mu}^2$ in $\mathcal{L}^1$. For completeness we prove it here. Let $\omega > 0$ be fixed, $V_l' = V_l \mathbf{1}_{V_l \leq \omega}$ and $V_l'' = V_l - V_l'$. By (4.2),

$$\limsup_{n\to\infty} \left\| \sum_{l=0}^{\infty} c_{n,l}^2 (V_l' - \mathbf{E} V_l') \right\|^2 \leq \limsup_{n\to\infty} \sum_{l,l'=0}^{\infty} c_{n,l}^2 c_{n,l'}^2 |\mathrm{cov}(V_0', V_{l-l'}')| = 0$$

since $\lim_{k\to\infty} \mathrm{cov}(V_0', V_k') = 0$. Therefore, using $V_l = V_l' + V_l''$, again by (4.2),

$$\limsup_{n\to\infty} \mathbf{E} \left| \sum_{l=0}^{\infty} c_{n,l}^2 (V_l - \mathbf{E} V_l) \right| \leq \limsup_{n\to\infty} \mathbf{E} \left| \sum_{l=0}^{\infty} c_{n,l}^2 (V_l'' - \mathbf{E} V_l'') \right|$$

$$\leq 2 \limsup_{n\to\infty} \sum_{l=0}^{\infty} c_{n,l}^2 \mathbf{E} V_l'' \to 0 \text{ as } \omega \to \infty.$$

Under (4.2), for any $\delta > 0$, $\sum_{l=0}^{\infty} \mathbf{E}\{|c_{n,l}^2 \eta_{m_2-l}^2|\mathbf{1}_{|c_{n,l}\eta_{m_2-l}|\geq \delta}\} \to 0$. So the finite dimensional convergence holds. By Proposition 4 of Dedecker and Doukhan [8],

$$\|S_n\|_q^2 \leq 2q \|\eta_0\|_q^2 \sum_{j=1}^{\infty} (A_j - A_{j-n})^2 = O(\sigma_n^2).$$

By Theorem 2.1 of Taqqu [32], the tightness follows. (2°) Note that $\|S_n\| \sim \|\eta_0\| \sigma_n$, the conclusion similarly follows.

[b] The finite dimensional convergence follows in the same manner as [a]. For the tightness, let $1 \leq m_1 < m_2 \leq n$, by Proposition 4 in [8],

$$\|S_{m_2}^\diamond - S_{m_1}^\diamond\|_q^2 \leq 2q \|\eta_0\|_q^2 \sum_{j=0}^{m_2-1} (A_j - A_{j-(m_2-m_1)})^2 = O(\tilde\sigma_{m_2-m_1}^2).$$

With the above inequality, using the same argument as in Theorem 2.1 of Taqqu [32], we have for any $0 \leq t_1 \leq t \leq t_2 \leq 1$, there exists a generic constant $C$ (independent of $n, t_1, t$ and $t_2$), such that for $\beta \in (1/2, 1)$,

$$\mathbf{E}(|S_{\lfloor nt \rfloor}^\diamond - S_{\lfloor nt_1 \rfloor}^\diamond| |S_{\lfloor nt_2 \rfloor}^\diamond - S_{\lfloor nt \rfloor}^\diamond|) \leq C \tilde\sigma_n^2 (t_2 - t_1)^{3-2\beta}$$

and for $\beta \in (1, 3/2)$,

$$\mathbf{E}|S_{\lfloor nt \rfloor}^\diamond - S_{\lfloor nt_1 \rfloor}^\diamond|^{q/2} |S_{\lfloor nt_2 \rfloor}^\diamond - S_{\lfloor nt \rfloor}^\diamond|^{q/2} \leq C \tilde\sigma_n^q (t_2 - t_1)^{q(3/2-\beta)}.$$

Thus the tightness follows from Theorem 15.6 in [3]. □

**Remark 4.1.** Under (1°) of Lemma 4.1, the moment condition $\eta_j \in \mathcal{L}^q$, $q > (3/2 - \beta)^{-1}$, is optimal. and it can not be reduced to $\eta_j \in \mathcal{L}^{q_0}$, $q_0 = (3/2 - \beta)^{-1}$. Consider the case in which $\eta_i$ are iid symmetric random variables and $\mathbf{P}(|\eta_0|^{q_0} \geq g) \sim g^{-1}(\log g)^{-2}$ as $g \to \infty$. Then $\eta_j \in \mathcal{L}^{q_0}$. Let $\ell(n) = 1/\log n$, $n > 3$. Elementary calculations show that $\sigma_n^{-1} \max_{1\leq j\leq n} |\eta_j| \to \infty$ in probability. Let $Y_j' = Y_j + \eta_j - \eta_{j-1}$ and $S_i' = \sum_{j=1}^{i} Y_j'$. Then the coefficients $a_j'$ of $Y_j'$ also satisfy the conditions in Lemma 4.1. The two processes $\sigma_n^{-1} S_{\lfloor nt \rfloor}$ and $\sigma_n^{-1} S_{\lfloor nt \rfloor}'$, $0 \leq t \leq 1$, cannot both converge weakly to fractional Brownian motions. If so, since $\max_{j\leq n} |\eta_j - \eta_0| \leq \max_{j\leq n} |S_j| + \max_{j\leq n} |S_j'|$, we have $\max_{j\leq n} |\eta_j| = O_{\mathbf{P}}(\sigma_n)$, contradicting $\sigma_n^{-1} \max_{1\leq j\leq n} |\eta_j| \to \infty$ in probability. Similar examples are given in [41, 43].

*Proof of Theorem 2.1.* Let $a_j = \Gamma(j+d)/\{\Gamma(d)\Gamma(j+1)\}$, $j \geq 0$, and $A_k = \sum_{i=0}^{k} a_i$ if $k \geq 0$ and $0$ if $k < 0$. Note that $\beta = 1 - d$. Then $X_t = \sum_{j=0}^{\infty} a_j u_{t-j}$. By (2.1), $\zeta_j = \sum_{i=j}^{\infty} \mathcal{P}_j u_i \in \mathcal{L}^q$. Let $M_n = \sum_{i=1}^{n} \zeta_i$, $S_n = \sum_{j=1}^{n} Y_j$, $Y_j = \sum_{i=0}^{\infty} a_i \zeta_{j-i}$, $U_n = \sum_{i=1}^{n} u_i$ and $R_n = T_n - S_n$. By Theorem 1 in [39], $\|U_n - M_n\|_q = o(\sqrt{n})$. By Karamata's theorem and summation by parts, we have

$$\|R_m\|_q \leq \left\| \sum_{i=0}^{3m} (A_i - A_{i-m})(u_{m-i} - \zeta_{m-i}) \right\|_q$$
$$+ \left\| \sum_{i=3m+1}^{\infty} (A_i - A_{i-m})(u_{m-i} - \zeta_{m-i}) \right\|_q$$
$$= \sum_{i=1}^{3m} |(A_i - A_{i-m}) - (A_{i-1} - A_{i-1-m})|o(\sqrt{i})$$
$$+ \sum_{i=3m+1}^{\infty} |(A_i - A_{i-m}) - (A_{i-1} - A_{i-1-m})|o(\sqrt{i}) = o(\sigma_m).$$

By Proposition 1 in [39],

$$\left\| \max_{m \leq 2^k} |R_m| \right\|_q \leq \sum_{r=0}^{k} 2^{(k-r)/q} \|R_{2^r}\|_q = \sum_{r=0}^{k} 2^{(k-r)/q} o(\sigma_{2^r}) = o(\sigma_{2^k}),$$

since $q > 2/(2d+1)$. So the limit of $\{T_{\lfloor nt \rfloor}/\sigma_n, 0 \leq t \leq 1\}$, if exists, is equal to the limit of $\{S_{\lfloor nt \rfloor}/\sigma_n, 0 \leq t \leq 1\}$. By Lemma 4.1, the latter has a weak limit. So (2.2) follows. □

*Proof of Theorem 2.2.* As in the proof of Theorem 2.1, let $S_i^\diamond = \sum_{j=1}^{i} Y_j^\diamond$ and $R_m^\diamond = T_m^\diamond - S_m^\diamond$. By Karamata's theorem,

$$\|R_{m+l}^\diamond - R_l^\diamond\|_q \leq \sum_{j=1}^{m} |A_j - A_{j-1}|o(\sqrt{j}) = o(\tilde{\sigma}_m).$$

Again by the maximal inequality (Proposition 1 in [39]),

$$\left\| \max_{m \leq 2^k} |R_m^\diamond| \right\|_q \leq \sum_{r=0}^{k} \left[ \sum_{j=1}^{2^{k-r}} \|R_{2^r j}^\diamond - R_{2^r(j-1)}^\diamond\|_q^q \right]^{1/q}$$
$$= \sum_{r=0}^{k} 2^{(k-r)/q} o(\tilde{\sigma}_{2^r}) = o(\tilde{\sigma}_{2^k}),$$

which proves the theorem in view of Lemma 4.1. □

*Proof of Theorem 3.1.* If (3.3) holds, by the continuous mapping theorem, Theorem 2.1 entails (3.4) and (3.5). In the sequel we shall prove (3.3). Note that

$$w_{n,l}^2 = \frac{1}{n} \sum_{j=1}^{n} X_j^2 + \frac{2}{n} \sum_{j=1}^{l} \left(1 - \frac{j}{l+1}\right) \sum_{i=1}^{n-j} X_i X_{i+j} - \bar{X}_n^2$$
$$- \frac{2\bar{X}_n}{n} \sum_{j=1}^{l} \left(1 - \frac{j}{l+1}\right) \sum_{i=1}^{n-j} (X_i + X_{i+j}) + \frac{2\bar{X}_n^2}{n} \sum_{j=1}^{l} \left(1 - \frac{j}{l+1}\right) (n-j)$$
$$=: I_{1n} - \bar{X}_n^2 + I_{2n}.$$

Since $|\bar{X}_n| = O_{\mathbf{P}}(n^{d-1/2})$, $l^{-2d}(\bar{X}_n^2 + |I_{2n}|) = O_{\mathbf{P}}((l/n)^{1-2d}) = o_{\mathbf{P}}(1)$. Thus it suffices to show that $l^{-2d}I_{1n} \to_{\mathbf{P}} \kappa_1^2(d)$. Let $V_j = \sum_{i=1}^j X_i$, $\tilde{V}_j = \sum_{i=n-j+1}^n X_i$, $1 \leq j \leq l$, then a straightforward calculation shows that $I_{1n} = J_{1n} + J_{2n}$, where

$$J_{1n} := \frac{1}{(l+1)n} \sum_{i=l+1}^n \left(\sum_{j=i-l}^i X_j\right)^2, \quad J_{2n} := \frac{1}{n(l+1)} \sum_{j=1}^l (V_j^2 + \tilde{V}_j^2).$$

Corollary 2.1 implies that $\sum_{j=1}^l V_j^2 = O_{\mathbf{P}}(l^{2d+2})$. Since $\sum_{j=1}^l \tilde{V}_j^2$ has the same distribution as $\sum_{j=1}^l V_j^2$, we have $J_{2n} = O_{\mathbf{P}}(l^{2d+2}/(ln)) = o_{\mathbf{P}}(l^{2d})$.

It remains to show $l^{-2d}J_{1n} \to_{\mathbf{P}} \kappa_1^2(d)$. Let $\zeta_j = \sum_{i=j}^\infty \mathcal{P}_j u_i \in \mathcal{L}^q$, $M_n = \sum_{i=1}^n \zeta_i$, $U_n = \sum_{i=1}^n u_i$ and $r_n = \sup_{j \geq n} \|U_j - M_j\|/\sqrt{j}$. Then $r_n \to 0$ and it is non-increasing. Let $L = \lfloor \min\{\sqrt{nl}, l(r_{\sqrt{l}})^{1/(-1-2d)}\} \rfloor$. Then $l = o(L)$ and $L = o(n)$. Let

$$W_{j,l} = \sum_{i=0}^L (A_i - A_{i-l-1}) u_{j-i}, \quad Q_{j,l} = \sum_{i=0}^L (A_i - A_{i-l-1}) \zeta_{j-i}$$

and $S_{j,l} = \sum_{i=j-l}^j X_i$. Then $J_{1n} = \frac{1}{n(l+1)} \sum_{j=l+1}^n S_{j,l}^2$. Since $l/L \to 0$,

(4.3)
$$\|S_{j,l} - W_{j,l}\| = \left\|\sum_{i=L+1}^\infty (A_i - A_{i-l-1}) u_{j-i}\right\|$$
$$\leq \left[\sum_{i=L+1}^\infty (A_i - A_{i-l-1})^2\right]^{1/2} \sum_{t=0}^\infty \|\mathcal{P}_0 u_t\| = o(\sigma_l).$$

By the definition of $L$, using summation by parts, we have

(4.4)
$$\|W_{j,l} - Q_{j,l}\| \leq \sum_{i=1}^L |(A_i - A_{i-l-1}) - (A_{i-1} - A_{i-l-2})| r_i \sqrt{i}$$
$$\leq r_{\sqrt{l}} \sum_{i=1+\sqrt{l}}^L 2|a_i| \sqrt{i} + \sum_{i=1}^{\sqrt{l}} 2|a_i| r_i \sqrt{i} = o(\sigma_l).$$

Now we shall show

(4.5)
$$\mathbf{E}\left|\sum_{j=1}^n \{Q_{j,l}^2 - \mathbf{E}(Q_{j,l}^2)\}\right| = o(n\sigma_l^2).$$

Since $\tau_m := \mathbf{E}|\mathbf{E}(\zeta_0^2|\mathcal{F}_{-m}) - \mathbf{E}(\zeta_0^2)| \to 0$ as $m \to \infty$,

(4.6) $\quad \mathbf{E}|\mathbf{E}(Q_{2Lk,l}^2|\mathcal{F}_{2Lk-2L}) - \mathbf{E}(Q_{2Lk,l}^2)| \leq \sum_{i=0}^L (A_i - A_{i-l-1})^2 \tau_{2L-i} = o(\sigma_l^2).$

Let $D_k = Q_{2Lk,l}^2 - \mathbf{E}(Q_{2Lk,l}^2|\mathcal{F}_{2Lk-2L})$, $k = 1, \ldots, b = \lfloor n/(2L) \rfloor$. Set $C_q = 18q^{3/2}(q-1)^{-1/2}$ and $q' = \min(q, 4)$. By Burkholder's inequality,

(4.7)
$$\left\|\sum_{k=1}^b D_k\right\|_{q/2} \leq C_{q/2} b^{2/q'} \|D_k\|_{q/2} \leq 2C_{q/2} b^{2/q'} \|Q_{2Lk,l}\|_q^2$$
$$= b^{2/q'} O(\sigma_l^2) = o(b\sigma_l^2).$$

Thus (4.5) follows from (4.6) and (4.7). By (4.3), (4.4) and (4.5),

$$\mathbf{E}\left|\sum_{j=1}^{n}\{S_{j,l}^2 - \mathbf{E}(Q_{j,l}^2)\}\right| \leq \sum_{j=1}^{n} \mathbf{E}|S_{j,l}^2 - Q_{j,l}^2| + o(n\sigma_l^2) = o(n\sigma_l^2)$$

which completes the proof since $\mathbf{E}(Q_{j,l}^2) \sim \sigma_l^2 \|\zeta_0\|^2$.  □